# A Generalized Analytical Heat Transfer Model for Enhanced Geothermal Systems: Capturing Fracture Interactions and Correcting Classical Optimistic Predictions

Nelson Barros-Galvis[1], Christine Ehlig-Economides[2], Cristi Darley Guevara[3]

**Abstract**

Numerical analytical heat transfer models play a critical role in geothermal design and feasibility studies. Classical solutions, such as Gringarten et al. (1975), often assume simplified conditions and systematically overestimate thermal performance, potentially leading to unrealistic engineering decisions.

This study presents a generalized analytical model for enhanced geothermal systems (EGS) that captures the interactions between fractures while preserving analytical tractability. The formulation is based on Green's functions and reproduces realistic thermal behavior under conditions representative of fractured reservoirs. Notably, the solution is sufficiently simple and computationally efficient that it can be implemented in a standard spreadsheet without requiring Laplace-space transformations numerical algorithms, making it directly applicable for engineering practice.

The model is validated against numerical simulations using CMG STARS and Volsug software, showing close agreement for temperature evolution, including the effects of interacting fractures. Compared to classical analytical approaches, the proposed model corrects optimistic biases and provides more reliable predictions for production temperature and energy recovery.

The results have direct implications for geothermal feasibility studies, well design, and power forecasting**,** bridging the gap between legacy analytical models and numerical/commercial engineering tools. This approach enables engineers and decision-makers to plan geothermal projects with greater confidence in reservoir performance predictions, particularly in enhanced and fractured geothermal systems.

Building on the analytical framework originally introduced by Gringarten et al. (1975).  The proposed formulation generalizes classical heat transfer solutions to account for interacting fractures while retaining analytical tractability.

**Keywords:** enhanced geothermal system; heat transfer; analytical model; Gringarten's model,

---

[1] Corresponding Author:  University of Houston, 5000 Gulf Freeway, Building 9, Room 141, Houston, USA *nelbarrosgalvis@gmail.com , nebarros@central.uh.edu*

[2] University of Houston, 5000 Gulf Freeway, Building 9, Room 147, Houston, USA *ceconomi@central.uh.edu*

[3] Arizona State University, ASU-Helios Decision Center, University Dr & Mill Av, Tempe Arizona, USA *cristi.guevara@asu.edu*

**1. Introduction**

Enhanced Geothermal Systems (EGS) are a promising sustainable source of renewable energy, capable of providing reliable electricity from deep geothermal resources. Prediction of reservoir thermal behavior and power is critical for feasibility studies, well design, and long-term energy planning. Analytical models are particularly valuable in early-stage assessments due to their simplicity. computational efficiency, and tractability, allowing rapid evaluation of multiple scenarios.

The classical analytical solution proposed by Gringarten et al. (1975) remains widely cited and provides a framework to estimate heat transfer in fractured hot dry rock. However, it is known to produce optimistic predictions, potentially overestimating reservoir performance. Practical implementation of the Gringarten model often requires specialized numerical algorith such as Stephest to implement the analytical solution, which can limit its accessibility in routine engineering evaluations.

Numerical approaches have been developed to provide more realistic estimates, and to overcome of classical analytical solutions. Zeilani et al. (2021) employed CMG simulations to model EGS using triplet-wells, offering more realistic predictions of rock behavior, including water outlet temperature evolution. Interestingly, the CMG model tends to be slightly pessimistic, highlighting potential underperformance in design calculations compared to classical analytical estimates. Similarly, Stacey, (2025) implemented a numerical model using Volsug software, comparing its predictions against Gringarten's analytical solution. Volsug provided validation for model behavior but was not benchmarked against field data. While highly accurate, both CMG STARS and Volsug approaches are computationally intensive and require specialized commercial software, restricting their direct use in early-stage or spreadsheet-based engineering studies.

In this context, there is a clear need for an analytical model that combines the computational simplicity of classical solutions with the realistic behavior captured by numerical simulations. This study presents a generalized analytical heat transfer model for EGS that captures the interactions between fractures, corrects the optimistic bias of classical formulations, and preserves analytical tractability. The solution is sufficiently simple to be implemented directly in standard spreadsheets, making it immediately applicable for engineering practice.

Validation of the proposed model is performed against Gringarten's analytical predictions, CMG simulations (Zeilani et al.)**,** and Volsug numerical results (Stacey). This comprehensive comparison demonstrates that the model provides accurate and practical thermal predictions, bridging the gap between classical analytical and numerical approaches. The model therefore offers engineers and decision-makers a reliable tool for forecasting production temperatures, thermal performance, and overall geothermal feasibility.

The present work builds upon this analytical tradition by extending the classical framework to capture thermal interactions between fractures and finite heat transfer effects. By retaining an analytical formulation grounded in Green's functions, the proposed model preserves the interpretability and computational efficiency of classical solutions while addressing key sources of optimistic bias identified in earlier formulations.

## 2. Analytical heat transfer models in geothermal engineering

Analytical heat transfer models have long played an important role in geothermal engineering due to their conceptual clarity, computational efficiency, and ability to provide rapid estimates of reservoir thermal behavior. Such models are particularly valuable during the early stages of geothermal project evaluation, where multiple design scenarios must be screened efficiently before committing to detailed numerical simulations. In the context of fractured and enhanced geothermal systems, analytical formulations offer physical insight into the coupled effects of fluid flow and heat transfer, while enabling parametric studies that are often impractical with fully numerical approaches.

### 2.1 Classical Analytical Models

The formulation proposed by Gringarten et al. (1975) represents a cornerstone in geothermal reservoir engineering. The model is derived using Laplace transforms, leading to an exact analytical solution in Laplace space for heat extraction from fractured hot dry rock systems. Due to the mathematical complexity of the resulting expressions, practical evaluation of the solution in the time domain requires numerical inversion of the Laplace transform, commonly performed using the method of Papoulis. Despite this complexity, the model has been widely adopted as a reference for feasibility studies, conceptual design, and benchmarking of numerical simulators in geothermal engineering.

Furthermore, in the classical formulation of Gringarten et al. (1975), heat transfer is modeled under a set of boundary conditions and expressed using dimensionless variables. The reservoir is represented as a hot dry rock system containing an infinite series of parallel, vertical fractures that are equidistant and thermally isolated from one another, such that thermal interaction between adjacent fractures is neglected. The surrounding rock matrix is assumed to be laterally infinite, and heat transfer is governed solely by conductive heat exchange from the hot dry rock matrix to each fracture.

These assumptions inherently imply low fluid velocities and, consequently, effectively long residence times within the fracture system, conditions that favor near-complete thermal equilibration between the circulating fluid and the surrounding rock, resulting in a heat transfer regime that is strongly conduction-dominated. The initial temperature distribution is defined as a function of depth, $z$, through a constant geothermal gradient, which serves as the thermal boundary condition for the system. The thermal breakthrough approximately is at 18 years approxiamadely.

Within this framework, fluid flow and heat exchange are treated as uniform along the fracture planes, and the evolution of produced fluid temperature is described in terms of dimensionless groups that enable generalization of the solution across different fracture spacings and reservoir scales. The formulation is strictly valid for an infinite number of fractures, a limitation that is evident in the classical type curves of dimensionless water outlet temperature versus dimensionless time where the thermal breakthrough is observed at approximately 18 years. While these assumptions facilitate analytical tractability and the construction of type curves, they represent an idealized upper bound on thermal performance and do not capture fracture–fracture interference. or finite thermal boundaries characteristic of modern enhanced geothermal systems.

**2.2 Numerical Models as Reference**

To overcome the limitations of analytical models, numerical simulations have been widely used. Zeinali (2021) employed CMG STAR to simulate EGS resource, providing predictions that tend to be more conservative compared with classical analytical solutions. Stacey (2025) used Volsug software to compare its results with Gringarten's analytical findings, and it was noted that the thermal breakthrough occurs at approximately 7 years. The internal methods and governing equations of CMG STAR and Volsug are not publicly disclosed. In this work, the generalized analytical model is validated against CMG STAR, Volsug, and the original Gringarten results, showing close agreement across these numerical and analytical approaches. While numerical simulations are computationally intensive and require specialized software, the analytical formulation presented here allows direct evaluation and spreadsheet-based calculations for early-stage engineering studies.

**2.3 Heat Transfer in Fractured Reservoirs**

A key factor influencing thermal behavior in fractured hot dry rock and enhanced geothermal systems is the interaction between fractures. Heat transfer in one fracture is affected by the spatial distribution and connectivity of neighboring fractures. Classical analytical models, including Gringarten's solution, assume simplified fracture geometries and neglect thermal interaction between fractures. These simplifications can lead to overestimation of resource temperature and thermal performance, highlighting the need for generalized analytical models that explicitly account for fracture–fracture interactions.

**2.4 Motivation for a Generalized Analytical Model**

The limitations of classical analytical and numerical approaches highlight the need for a model that combines tractable analytical solutions with physically realistic behavior**.** Despite decades of development, there is a clear gap in existing literature for a model that is:

- Analytical, allowing rapid evaluation without heavy computational requirements
- Capable of capturing thermal interactions between fractures, which strongly influence thermal evolution

- Correcting the optimistic bias of classical formulations
- Directly implementable in standard engineering tools, such as spreadsheets

**3. Origin of optimistic bias in classical formulations**

The classical analytical solution proposed by Gringarten et al. (1975) provides a tractable framework to estimate the evolution of produced fluid temperature in fractured geothermal reservoirs. However, the formulation systematically leads to optimistic thermal predictions, which can be traced to several fundamental physical assumptions.

The model assumes isolated infinite vertical fractures, implicitly neglecting thermal interactions between adjacent fractures. In realistic enhanced geothermal systems, heat extraction occurs through networks of hydraulically and thermally connected fractures, where interference effects significantly influence temperature evolution. By treating fractures as isolated entities, the classical formulation overestimates the amount of heat available to each fracture.

Additionally, the model is implicitly formulated under conditions of extremely low fluid velocities, which imply very long residence times within the fracture system. These conditions favor nearly complete thermal equilibration between the rock and the circulating fluid, further contributing to optimistic production temperature predictions.

A key source of bias lies in the behavior of the dimensionless parameter $\beta$, which is associated with the thermal boundary condition and the geothermal gradient. For the classical solution to remain valid, $\beta$ must approach zero, effectively implying unlimited heat transfer from the reservoir to the fluid. This assumption represents an idealized upper bound rather than a realistic operating condition for modern EGS, where finite heat transfer rates and non-negligible thermal resistance are present.

Finally, practical implementation of the analytical solution often requires Laplace-space transformations and numerical inversion techniques. While these approaches preserve analytical tractability, they do not modify the underlying physical assumptions responsible for optimistic thermal predictions.

While these assumptions preserve analytical simplicity and yield closed-form solutions, they result in temperature and performance predictions that are higher than those expected under realistic geothermal operating conditions. Recognizing the physical origin of this optimistic bias motivates the development of generalized analytical models capable of capturing fracture interactions and finite heat transfer effects while maintaining analytical tractability.

**4. Generalized analytical heat transfer model**

To address the limitations identified in classical analytical formulations and numerical simulations, this study introduces a generalized analytical heat transfer model for fractured hot dry rock and enhanced geothermal systems. The model is designed to retain the analytical tractability of classical solutions while incorporating physically realistic features that are essential for modern EGS, particularly thermal interactions between fractures.

Unlike classical formulations that assume thermally isolated fractures, the proposed model explicitly accounts for fracture–fracture thermal interference through the surrounding rock matrix. The formulation is derived under conduction-dominated heat transfer in the rock matrix and convective heat exchange within the fractures, consistent with typical operating conditions in granitic geothermal reservoirs. The resulting analytical solution provides realistic predictions of produced fluid temperature and thermal breakthrough while remaining directly implementable in standard engineering tools, such as spreadsheets.

**4.1 Physical framework and assumptions**

The geothermal system considered in this study represents an enhanced geothermal system (EGS) developed in granitic hot dry rock. Heat extraction occurs through a finite set of hydraulically stimulated fractures embedded within a finite reservoir domain. Heat transfer in the surrounding rock matrix is assumed to be conduction-dominated, while convective heat transport occurs within the fractures due to fluid circulation.

The fracture network is represented by a series of discrete, parallel, horizontal fractures of finite thickness, intersecting a horizontal well, and distributed along the length of the stimulated reservoir volume. Unlike classical analytical formulations that assume an infinite, periodically spaced array of thermally isolated fractures, the present framework explicitly considers a finite number of fractures within a finite reservoir length. As a result, thermal interactions between fractures naturally arise through the surrounding rock matrix and are accounted for in the formulation.

A key physical consideration in finite EGS reservoirs is the presence of thermal boundaries at the edges of the stimulated domain. Fractures located near these boundaries experience asymmetric heat drainage, as the available rock volume for heat extraction is reduced compared to that of interior fractures. This boundary effect leads to earlier thermal depletion and diminished long-term thermal contribution from edge fractures if uniform spacing is imposed in reservoir length.

Other methodology is, the reservoir is decomposed into a series of longitudinal blocks, each drained by a uniformly spaced fracture cluster. Within each block, classical analytical solutions apply, while inter-block thermal interaction is captured through effective boundary conditions and volume-weighted averaging.

To mitigate this effect, the framework allows adjustment of the fracture array extent relative to the reservoir boundaries, enabling uniformly spaced fractures to access comparable thermal drainage volumes. By repositioning fractures adjacent to the physical boundaries, the effective thermal drainage volume of edge fractures can be made comparable to that of interior fractures, enabling boundary fractures to behave thermally as intermediate fractures.

The reservoir is decomposed into a series of longitudinal blocks, each drained by a uniformly spaced fracture cluster. Within each block, classical analytical solutions apply, while inter-block thermal interaction is captured through effective boundary conditions and volume-weighted averaging. This design concept is naturally captured within the generalized analytical formulation and provides a physically consistent means to account for finite-domain effects in EGS design.

**5. Governing Equations and Green's Function Formulation**

Cartesian coordinates are adopted in this study to ensure methodological consistency and facilitate direct comparison with existing analytical and numerical models. The classical formulation of Gringarten et al. (1975), as well as the numerical implementations in CMG STARS and Volsung software, are all formulated in rectangular coordinate systems. By adopting the same coordinate framework, differences in thermal predictions can be directly attributed to physical assumptions and model formulation, rather than geometric or coordinate-system effects.

The following assumptions are made:

1. The parallel horizontal fractures have uniform aperture, and the rock matrix is impermeable.
2. The rock and fluid thermal conductivities, heat capacities, densities are constant.
3. The fracture spacing is equidistant.
4. The stable fluid transport is given in hydraulic fracture.
5. Multiple fractures generate heat interaction related to their spacing.
6. Thermal stresses and granitic rock creep have not been considered in this study.

**5.1 Governing equations**

Heat transfer in the enhanced geothermal system considered in this study is governed by conductive heat transport in the surrounding granitic rock matrix and convective heat transport within the hydraulically stimulated fractures. Under EGS operating conditions, heat transport in the rock matrix is conduction dominated and transient due to the thermal inertia of the rock mass. In contrast, heat transport within the fractures is dominated by fluid advection and is treated as quasi steady-state along the fracture plane, reflecting the comparatively short fluid residence times relative to the characteristic timescale of conductive heat diffusion in the surrounding rock.

Rock heat conduction:
The temperature evolution in the rock matrix is governed by the transient heat conduction equation

$$\frac{\partial T_r}{\partial t} = \alpha \frac{\partial^2 T_r}{\partial y^2}, \quad y > 0, \ t > 0 \qquad [1]$$

The temperatures must satisfy the following conditions:

Initial condition:

$$T_r(y, 0) = T_0 \qquad [2]$$

Dirichlet boundary condition for $y = 0$:

$$T_r(0, t) = T_f(x, t) \qquad [3]$$

Where

$T_0$ is initial rock temperature.

$T_r$ is the rock temperature

$\alpha$ is diffusivity constant, $\alpha = \frac{k_r}{\rho_r c_r}$

$\rho_r$ is the rock density

$c_r$ is the rock heat capacity

$T_f(x, t)$ is fluid temperature

$k$ is the rock thermal conductivity

## 5.2 Rock-fracture interfacial heat flux

Heat exchange between the rock matrix and the circulating fluid occurs at the rock–fracture interface and is governed by continuity of heat flux. Heat exchange between the circulating fluid and the surrounding granitic rock is described through an interfacial heat flux at the fracture–matrix boundary.

This formulation enables a natural coupling between the heat conduction in the rock and the steady convective heat transport within the fracture.

$$\rho_f c_f v b \frac{\partial T_f}{\partial x} = -q(x,t) \qquad [4]$$

$q(x,t)$ is the total interfacial heat flux

$c_f$ is the fluid heat capacity

$T_f$ is the fluid temperature

$v$ is the fluid velocity

$\rho_f$ is the fluid density

$b$ is the fracture aperture

### 5.3 Fourier's law

Under conduction-dominated conditions in the rock matrix, heat conduction rate through surface is given by Fourier's law.

$$q(x,t) = -k \left[ \frac{\partial T_r}{\partial y}\bigg|_{y=0} \right] \qquad [5]$$

$q(x,t)$ is the total interfacial heat flux.

$k$ is the rock thermal conductivity

$T_r$ is the rock temperature

### 5.4 Green's functions and mathematical model solution

The Green's function for a semi-infinite domain represents the fundamental solution of equations [1] and [4] subject to conditions [2] and [3], and the Fourier equation. By construction, the Green's function satisfies the heat equation everywhere in the domain and enforces the prescribed interfacial boundary condition at the fracture surface, while ensuring finite temperature perturbations.

The use of Green's functions allows the temperature field in the rock matrix to be expressed as a convolution of the interfacial heat flux history with the thermal response of the semi-infinite medium. This formulation naturally captures the transient diffusion of heat away from the fracture and provides a physically consistent framework for modeling heat exchange between the rock matrix and the circulating fluid.

The operator $\partial_t - \alpha_r \partial_y^2$ semi-infinite domain, $y > 0$, the Green's function $G(y, t; y', \tau)$ satisfies

$$\frac{\partial G}{\partial t} - \alpha \frac{\partial^2 G}{\partial y^2} = \delta(y - y')\delta(t - \tau) \tag{6}$$

For initial and boundary conditions:

$$G(0, t; y', \tau) = 0 \tag{7}$$
$$G(y, 0; y', \tau) = 0 \quad \text{para } t < \tau \tag{8}$$

The result of simultaneous equations solution [1] and [4] subject to conditions [2] and [3], and the Fourier equation is

$$T_f(x, t) = T_0 + (T_0 - T_{\text{inj}}) \text{erf}\left(\frac{k_r x}{2\rho_f c_f v b \sqrt{\alpha t}}\right) \tag{9}$$

The equation [9] describes the behavior of fluid outlet temperature in EGS resource for single fracture case.

Where

$T_f(x, t)$ is the fluid outlet temperature

$T_{\text{inj}}$ is the cold injection fluid temperature

For clarity, the governing equations are presented in the main text, while the detailed derivation of the analytical solution using Green's functions is provided in Appendix A.

### 5.4 Influence of the thermal diffusion radius on thermal interaction between hydraulic fractures

Influence thermal radius is the maximum distance from the fracture surface over which the thermal gradient is significantly affected during the operational lifetime of the project. It represents the thermal drainage zone of each fracture.

The transient temperature response induced by a continuous thermal source in a conductive medium can be expressed in closed form using the complementary error function. This represents the transient fundamental solution of the three-dimensional heat equation, derived using Green's functions, for a continuous point heat source of power P in a conductive medium.

$$\Delta T(r, t) = \frac{P}{4\pi k_r k} \text{erfc}\left(\frac{r}{2\sqrt{\alpha t}}\right) \tag{10}$$

This formulation naturally introduces a finite thermal radius that grows as, $r \propto 2\sqrt{\alpha t}$ enabling explicit representation of fracture–fracture thermal interaction.

The classical solution of the two-dimensional heat equation, derived using Green's functions, widely applied to fracture problems, planar heat sources, and line sources in conductive media is

$$\Delta T(x,y,t) = \frac{Q}{4\pi r}\left[E_1\left(\frac{r_1^2}{4\alpha t}\right) + E_1\left(\frac{r_2^2}{4\alpha t}\right)\right] \quad [11]$$

$E_1$ is exponential integral.

Characteristic radius is defined when exponent argument is 1.

$$\frac{r^2}{4\alpha t} = 1 \quad [12]$$

So that $r_{th} = r = \sqrt{4\alpha t}$

Where

$r_{th}$ is the influence thermal radius or characteristic distance

$t$ is characteristic time

$\alpha$ is the thermal diffusivity

Note that the argument of the error function in Equation 9 contains a similar expression describing the thermal radius and the characteristic time. For equidistant fractures, the onset of thermal interference occurs when the thermal radius equals half of the fracture spacing.

In an EGS with multiple equidistantly spaced fractures, the fluid outlet temperature is determined by the contribution of each fracture at a given time, accounting for the interaction effect between fractures characterized by the thermal radius of influence

## 6. Proposed model validation and comparison with Gringarten et al. (1975) and Stacey-Volsung models.

For the validation of the mathematical model, input data were taken from the Gringarten et al. (1975) analytical model for Valles Caldera, New Mexico, as show in Table 1. Notably, Gringarten's model was applied to cases with 1 and 10 fractures, considering fracture spacings of 40, 80, and 160. The geothermal system configuration follows a doublet pattern, consisting of one production well and one injection well. Table 1 shows the data used by Gringarten et al. (1975).

Table 1. Gringarten's model input for Valles Caldera New Mexico

| THERMAL PARAMETER | Value | Units (SI) | Value | Other Units |
|---|---|---|---|---|
| Initial rock temperatura, Tro | 300 | °C | 300 | °C |
| Inlet water temperature, Two | 65 | °C | 65 | °C |
| Water density, ρw | 1000 | kg/m³ | 1 | g/cm³ |
| Water specific heat, cw | 4184 | J/kg °C | 1 | cal/g °C |
| Rock thermal conductivity, Kr | 2.59408 | W/m °C | 0.0062 | cal/cm s °C |
| Rock density, ρr | 2650 | kg/m³ | 2.65 | g/cm³ |
| Rock specific heat, cr | 1046 | J/kg °C | 0.25 | cal/g °C |
| Fracture aperture, b | 0.00127 | m | 0.05 | in |
| Fracture height, H | 999.74 | m | 3280 | ft |
| Fracture lenght, 2Xf | 999.74 | m | 3280 | ft |
| Total volumetric flow rate, Q | 0.144 | m³/s | 78294.0 | bpd |
| Volumetric flow rate per fracture, Qf | 0.0144 | m³/s | 7829.4 | bpd |

First, Gringarten's model was reproduced considering a 50-year time to show technical details. Normally, geothermal project lifetime is considered between 15 – 30 years. The implementation involved analyzing the model equations and utilizing the Gaver-Stehfest algorithm for numerical Laplace inversion. Accordingly, we selected fracture spacings of 40 and 80 m for this analysis. Additionally, data obtained by Stacey (2025) using the Volsung simulator, based on 2D and 3D numerical modeling for fracture spacing of 40 m, were included. The resulting data are presented and compared against our analytical model in the Figure 1.

Figure 1 displays eight curves categorized by four colors. The solid orange line represents the outlet water temperature behavior for a single fracture as determined by our analytical model, while the dashed orange line shows the results from Gringarten et al. (1975). Both analytical models exhibit similar trends and concavity; the onset of temperature decline occurs at 0.32 and 0.30 years, respectively, with stabilization occurring around the 30-year mark. However, the final temperature at 50 years is 80°C for our model and 94°C for Gringarten's. This implies that extracting energy from a geothermal source using a single fracture is inefficient due to rapid thermal drawdown.

Figure 1 also illustrates the temperature behavior for a 40 m spacing, represented in gray. In this set, the solid line denotes our analytical model, the dashed line corresponds to Gringarten's model, and the circular markers represent Stacey's (2025) 2D numerical results. It should be noted that Stacey's 2D curve was digitized from the original publication; this curve aligns closely with Gringarten's model when implemented via the Gaver-Stehfest algorithm.

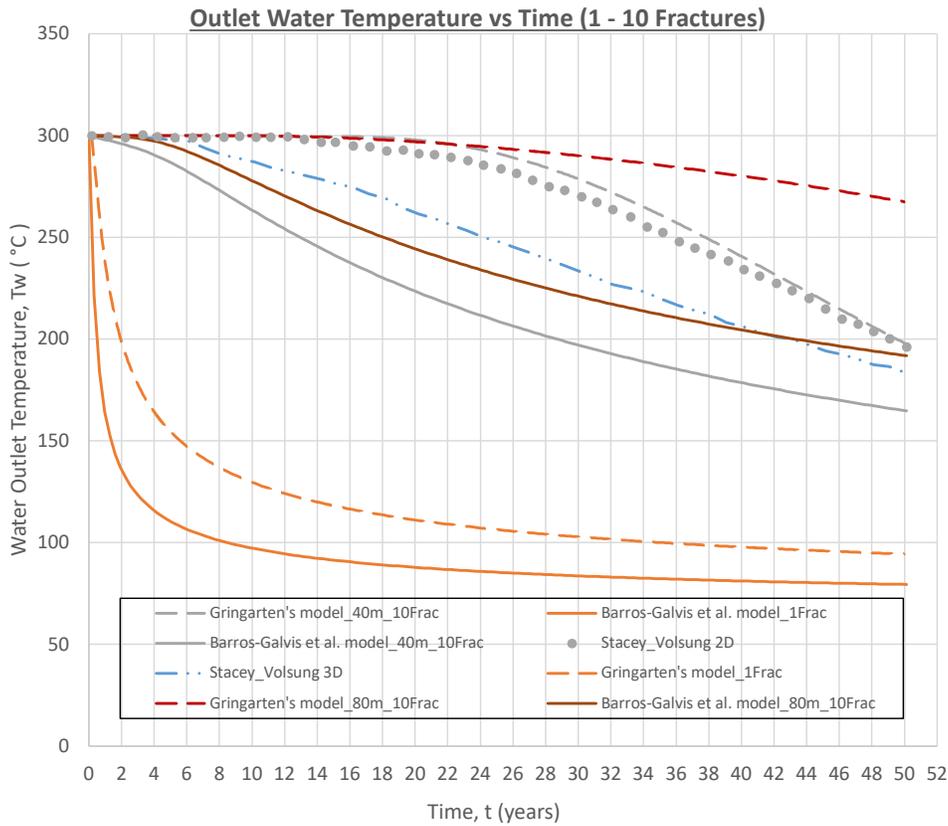

**Figure 1. Behavior of outlet water temperature vs time for 1 and 10 fractures using spacings of 40 and 80 m.**

Stacey reports that their 2D numerical model corresponds to Gringarten's analytical solution. In this case, the temperature decline begins at 18.1 years, reaching a final temperature of approximately 198°C at the 50-year mark. It should be noted that Gringarten's model assumes thermally isolated fractures with no interaction between them. In contrast, our analytical model, which accounts for thermal interference between fractures, shows that the decline starts much earlier at 2.6 years, with a fluid outlet temperature of approximately 165.3°C at a 50-year horizon.

Next, we analyze the 80 m fracture spacing, referenced in red. Our model predicts the onset of temperature decline at 4.2 years, with a final temperature of 192.4°C after 50 years. In contrast, Gringarten's model shows a temperature variation occurring at approximately 19 years, with a final fluid temperature of around 267°C. These discrepancies are significant; Gringarten's model presents an optimistic scenario where the geothermal reservoir energy remains virtually undepleted, suggesting that larger fracture spacing leads to higher outlet temperatures without accounting for the thermal radius of influence between fractures.

Special attention should be given to the 3D Stacey-Volsung case, shown in blue. Although this case uses a 40 m spacing and exhibits a graphical trend similar to our 40 m model—with a decline starting at year 7 and a final temperature of 180°C—its numerical results actually align more closely with our 80 m analytical model that includes thermal interaction. For both spacing configurations, our model provides a more conservative estimate regarding fluid outlet temperatures.

Thermal fracture interaction is detailed in Table 2, which presents calculated values for the thermal radius and the time to interference. For instance, with a fracture spacing of 80 m, a thermal front would take 54.214 years to traverse that distance. However, in the presence of an equidistant adjacent fracture, the neighboring thermal front propagates at the same rate, eventually meeting at a specific point in time. We define this as the 'interference time,' which corresponds to a thermal radius of interference that, in our mathematical model, equals half of the total thermal radius ($r_{thi} = r_{th}/2$). Physically, this interference zone behaves as a 'pseudo-fracture' at that moment, leading to an abrupt thermal decline. From an engineering perspective, designing EGS systems with spacings greater than 80 m is impractical, as interference occurs at 27.107 years—well within the standard 15 to 30-year evaluation horizon for geothermal projects.

**Table 2. Thermal radius and its interference time**

| Thermal radius, $r_{th}$ (m) | Time, t (years) | Interference time, $t_i$ (years) | Thermal interference radius, $r_{thi}$ (m) |
|---|---|---|---|
| 10 | 0.847 | 0.424 | 5 |
| 20 | 3.388 | 1.694 | 10 |
| 30 | 7.624 | 3.812 | 15 |
| 40 | 13.553 | 6.777 | 20 |
| 50 | 21.177 | 10.589 | 25 |
| 60 | 30.495 | 15.248 | 30 |
| 70 | 41.507 | 20.754 | 35 |
| 80 | 54.214 | 27.107 | 40 |

The thermal interference radius concept explains the conservative water outlet temperatures and the thermal drawdown rates observed in our analytical model, as presented in Figure 1. It is evident that Gringarten's model fails to account for the thermal radius of influence from a practical standpoint. Furthermore, the three-dimensional Stacey-Volsung model with 40 m spacing yields temperatures comparable to our analytical model with an 80 m spacing. This alignment is consistent with the fact that for a thermal radius of 80 m, the corresponding thermal interference radius is exactly 40 m.

## 6.1 Proposed model validation and comparison with Gringarten et al. (1975) and Zeinali-CMG models.

For the validation of the proposed model, a second comparative case was conducted using Gringarten's analytical model and the CMG-STARS numerical model presented by Zeilani et al. In their study, Zeilani et al. presented real-world cases featuring different volumetric flow rates per fracture (40, 80, and 160 bpd), with a fracture radius of 150 ft (45.72 m) and a fracture spacing of 130 ft (39.62 m) for a EGS system of 10 fractures. These simulation parameters are summarized in Table 3. EGS system is wells triplet configuration, namely, a pattern of production and injection wells in which each injector is surrounded by two producers.

Table 3. Zeilani's CMG model data for different rate: 40, 80, and 160 bpd

| PARAMETERS | Value | Units (SI) | Value | Other Units |
|---|---|---|---|---|
| Initial rock temperatura, Tro | 300 | °C | 300 | °C |
| Inlet water temperature, Two | 65 | °C | 65 | °C |
| Water density, ρw | 1000 | kg/m³ | 1 | g/cm³ |
| Water specific heat, cw | 4184 | J/kg °C | 1 | cal/g °C |
| Rock thermal conductivity, Kr | 2.59408 | W/m °C | 0.0062 | cal/cm s °C |
| Rock density, ρr | 2650 | kg/m³ | 2.65 | g/cm³ |
| Rock specific heat, cr | 1046 | J/kg °C | 0.25 | cal/g °C |
| Fracture aperture, b | 0.00127 | m | 0.05 | in |
| Fracture height, H | 91.44 | m | 300 | ft |
| Fracture lenght, 2Xf | 91.44 | m | 300 | ft |
| Fracture spacing | 39.62 | m | 130 | ft |
| Fracture radius | 45.72 | m | 150 | ft |

Figure 2 illustrates the outlet water temperature behavior over time. The curves are categorized by flow rates: green curves represent 160 bpd, red curves denote 80 bpd, and blue curves indicate 40 bpd. For each case, the results of the proposed analytical model, the numerical model, and Gringarten's model are depicted by a solid line, a dashed line with circular markers, and a dashed line without markers, respectively.

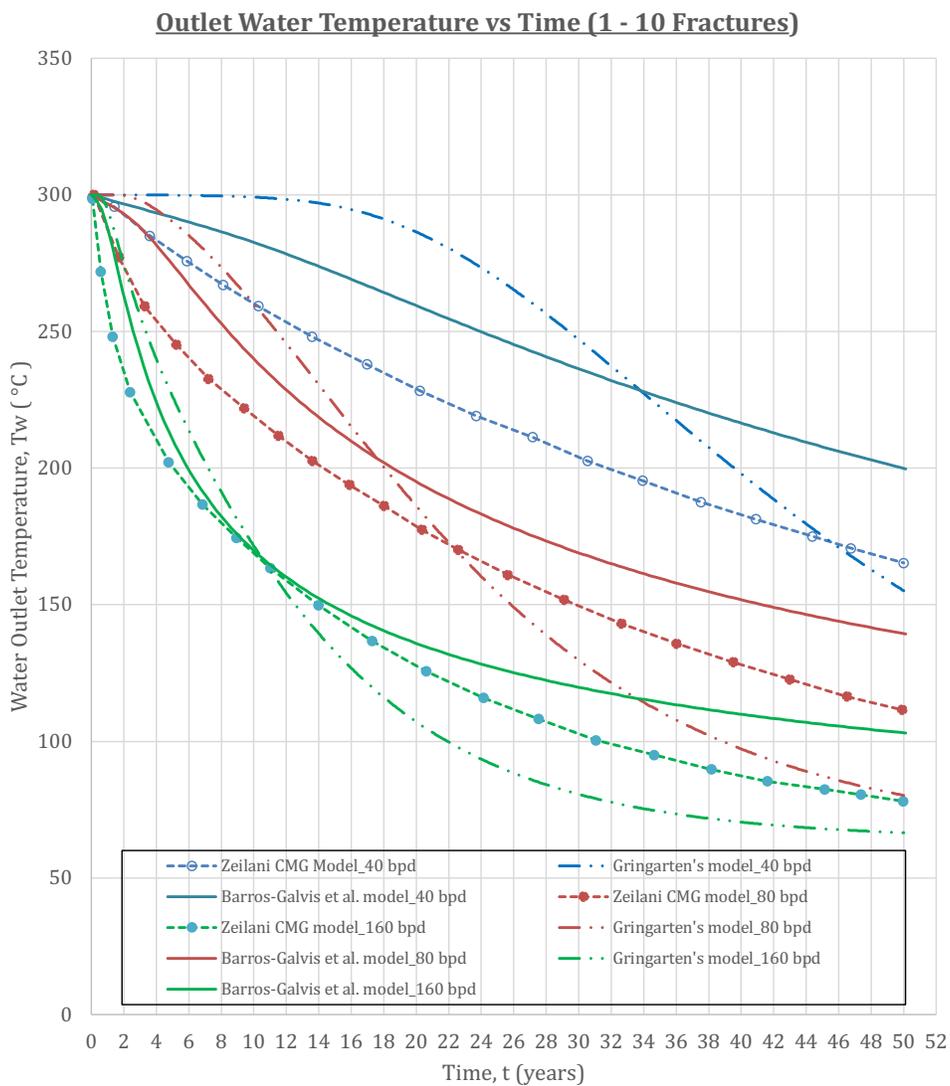

**Figure 2.** Behavior of outlet water temperature vs time for 10 fractures using spacings of 40 m and different flow rates (40, 80, and 160 bpd).

Initially, it is observed that across different fluid circulation rates, the outlet temperatures at 50 years are higher when the fluid velocity is lower. This is logically due to the increased residence time of the fluid within the fracture, which allows for better utilization of its heat capacity and enhanced heat transfer governed by convection-conduction mechanisms. Furthermore, the period during which the outlet temperature remains constant shortens as the volumetric flow rate increases. It is particularly noteworthy that both the proposed analytical model and the numerical model exhibit shorter periods of initial temperature stability. In Gringarten's model, this behavior is only approached by increasing the flow rate to 160 bpd. Moreover, our model proves to be conservative, distinguishing itself from Gringarten's by avoiding over-optimism, and from the CMG-STARS model by not being overly pessimistic.

The impact of volumetric flow on thermal drawdown is evident: higher rates significantly accelerate the onset of temperature decline. The final temperatures at the 50-year mark are summarized as follows:

- At 40 bpd: Gringarten (155°C), Proposed Model (199°C), CMG (165°C).
- At 80 bpd: Gringarten (80°C), Proposed Model (139°C), CMG (111°C).
- At 160 bpd: Gringarten (67°C), Proposed Model (103°C), CMG (78°C).

This comparison highlights that our analytical model serves as a balanced tool, offering conservative yet realistic forecasts for EGS performance across various operational scales, and avoiding the optimistic results of Gringarten and the more pessimistic outcomes of CMG.

Based on Table 2, the effect of the thermal radius of influence between equidistant fractures is clearly observed. In this case, with a fracture spacing of approximately 40 m, the corresponding thermal interference radius is 20 m. Interestingly, for a given spacing, increasing the fluid flow rate significantly shortens the time to the onset of thermal decline. This extreme fluid temperature depletion tends to mirror the behavior of a single, isolated fracture.

Ultimately, the Green's function-based model developed in this work bridges the gap between complex numerical modeling and practical engineering needs. It delivers a low-computational-cost alternative for EGS temperature analysis that is accessible, easy to implement in spreadsheets, and sufficiently accurate for long-term thermal assessment.

**8. Conclusions**

1. In conclusion, this study presents an analytical mathematical model developed using Green's functions that effectively determines the temporal evolution of fluid outlet temperatures in EGS. The model provides a robust alternative to high-cost numerical simulations, offering the significant advantage of practical implementation within standard spreadsheets without compromising physical accuracy.

2. This study demonstrates that the proposed analytical model accurately captures the thermal interaction between fractures, providing a more realistic representation of heat extraction compared to classical isolated-fracture models. The validation against the Stacey-Volsung 3D numerical simulator confirms that our analytical approach effectively approximates complex 3D behaviors by incorporating the thermal radius of influence.

3. While the Gringarten et al. (1975) model yields optimistic long-term temperatures by neglecting inter-fracture cooling, our model provides a conservative estimate essential for risk mitigation in EGS design. The significant discrepancy in outlet temperatures—up to 14°C in the single-fracture case and over 70°C in multi-fracture scenarios—underscores the risks of overestimating reservoir longevity.

4. Analysis of the interference time reveals that for an 80 m spacing, thermal fronts converge at 27.1 years. Given that geothermal projects typically operate on a 15 to 30-year horizon, designing systems with spacings greater than 80 m is technically inefficient, as the onset of rapid thermal drawdown occurs within the project's economic lifespan.

5. The identification of the interference zone as a "pseudo-fracture" provides a new physical framework for understanding abrupt temperature declines. This phenomenon marks the transition from independent heat extraction to a coupled regime where the rock matrix between fractures becomes thermally depleted.

6. The study confirms that for a given fracture spacing, thermal interference occurs at the midpoint. In the analyzed 40 m spacing case, the 20 m interference radius acts as a physical limit; once this threshold is reached at high flow rates, the multi-fracture system behaves thermally as a single, exhausted fracture, leading to rapid depletion.

7. Comparative analysis reveals that our model occupies a critical "middle ground" in geothermal forecasting. It corrects the optimistic bias of the Gringarten model (which ignores inter-fracture interference) without falling into the excessively pessimistic forecasts of some 3D numerical configurations. This positioning offers engineers a balanced and realistic tool for long-term production estimates.

8. The high degree of correlation between our proposed analytical model and the CMG-STARS numerical simulator results validates the model's robustness. By accurately replicating the thermal profiles across various flow rates (40, 80, and 160 bpd), the model proves to be a reliable alternative to computationally expensive numerical simulations for initial EGS assessments.

9. For a more robust engineering approach, further development of this model should include thermal stress and granite creep analysis, as these factors significantly influence the mechanical stability and permeability of EGS reservoirs

## Author Contributions

Nelson Barros-Galvis. The first author conceived, developed, solved, and validated the physical–mathematical model, performed comparisons with both analytical and numerical models, and led the writing of the manuscript.

Christine Ehlig-Economides. The second author provided overall guidance and supervision, contributed to the technical analysis and discussions, offered suggestions on the manuscript writing and figures.

Cristi Guevara. The third author reviewed and analyzed the physical–mathematical model solution.

# Appendix A

**Heat Conduction in the Rock Matrix**

The temporal evolution of temperature in the rock matrix is governed by the one-dimensional transient heat conduction equation

$$\frac{\partial T_r}{\partial t} = \alpha \frac{\partial^2 T_r}{\partial y^2}, \qquad y > 0, \ t > 0 \qquad [A1]$$

where $T_r$ is the rock temperature, $t$ is the time, $\alpha$ is diffusion constant.

The governing equation is subject to the following conditions.

Initial condition: $T_r(y, 0) = T_0,$ where $T_0$ is initial temperature of the rock.

Boundary condition at the fracture interface: $T_r(0,t) = T_f(x,t)$, where $T_f(x,t)$ is the fluid temperature within the fracture.

**Rock–Fracture Interfacial Heat Flux**

Heat exchange between the circulating fluid and the surrounding rock occurs at the rock–fracture interface and is governed by continuity of heat flux. The interaction between conductive heat transfer in the rock matrix and convective heat transport within the fracture is described through an interfacial heat flux at the fracture–matrix boundary. This formulation provides a natural coupling between the two domains.

$$\rho_f c_f v b \frac{\partial T_f}{\partial x} = -q(x,t) \qquad [A2]$$

where $q(x,t)$ is the total interfacial heat flux, $c_f$ is the fluid heat capacity, $T_f$ is the fluid temperature, $v$ is the fluid velocity, $\rho_f$ is the fluid density and $b$ is the fracture aperture

> **Commented [CDG1]:** heat capacity

**Fourier's Law**

Under conduction-dominated conditions in the rock matrix, the interfacial heat flux is described by Fourier's law,

$$q(x,t) = -k \left[ \frac{\partial T_r}{\partial y} \bigg|_{y=0} \right] \qquad [A3]$$

where $k$ is the thermal conductivity of the rock.

**Green's Function Formulation**

For a semi-infinite domain (y>0) with a Dirichlet boundary condition at y=0, the Green's function is constructed using the method of images,

$$G(y,t;\ y',\tau) = \frac{1}{\sqrt{4\pi\alpha(t-\tau)}} \left[ \exp\left(-\frac{(y-y')^2}{4\alpha(t-\tau)}\right) - \exp\left(-\frac{(y+y')^2}{4\alpha(t-\tau)}\right) \right] \qquad [A4]$$

The Green's function satisfies $\frac{\partial G}{\partial t} - \alpha \frac{\partial^2 G}{\partial y^2} = \delta(y-y')\,\delta(t-\tau)$

with the boundary condition $G(0,t;\ y',\tau) = 0$, and initial condition $G(y,0;\ y',\tau) = 0,\ for\ t < \tau$.

**Integral Solution for the Rock Temperature**

Using Green's integral representation, the solution of the heat conduction equation is

$$T_r(y,t) = \int_0^\infty G(y,t;y',0) T_r(y',0) dy' + \alpha \int_0^t \left[T_r(0,\tau) \frac{\partial G}{\partial y'}(y,t;0,\tau)\right]\bigg|_{y'=0} d\tau \quad \text{[A5]}$$

This expression is conveniently decomposed as

$$T_r(y,t) = T_{r1}(y,t) + T_{r2}(y,t) \quad \text{[A6]}$$

where $T_{r1}$ represents the contribution from the initial condition and $T_{r2}$ accounts for the boundary condition.

**Contribution of the Initial Condition**

The contribution due to the initial temperature field is

$$T_{r1}(y,t) = \int_0^\infty G(y,t;y',0) T_0 \, dy' \quad \text{[A7]}$$

Evaluating the Gaussian integrals yields

$$T_{r1}(y,t) = T_0 \, \text{erf}\left(\frac{y}{2\sqrt{\alpha t}}\right) \quad \text{[A8]}$$

This solution satisfies the initial condition $T_{r1}(y,0) = T_0$ for $y>0$ and vanishes at the boundary $y=0$, consistent with the separation of contributions.

**Physical interpretation**

The error function describes the diffusion of the initial rock temperature toward the fracture interface. For small times or large distances from the fracture, $\text{erf}(z) \to 1$, and the rock temperature remains close to $T_0$. Conversely, for large times or near the fracture, $\text{erf}(z) \to 0$ indicating progressive cooling of the rock adjacent to the fracture.

**Contribution of the Boundary Condition**

The contribution associated with the imposed boundary temperature is

$$T_{r2}(y,t) = \alpha \int_0^t T_f(x,\tau) \frac{\partial G}{\partial y'}\bigg|_{y'=0} d\tau \quad \text{[A9]}$$

Evaluating the derivative of the Green's function at $y' = 0$ gives

$$T_{r2}(y,t) = \frac{y}{\sqrt{4\pi\alpha}} \int_0^t \frac{T_f(x,\tau)}{(t-\tau)^{3/2}} \exp\left(-\frac{y^2}{4\alpha(t-\tau)}\right) d\tau \qquad [A10]$$

**Complete Solution for the Rock Temperature**

Combining both contributions, the rock temperature field is given by

$$T_r(y,t) = T_0 \operatorname{erf}\left(\frac{y}{2\sqrt{\alpha t}}\right) + \frac{y}{\sqrt{4\pi\alpha}} \int_0^t \frac{T_f(x,\tau)}{(t-\tau)^{3/2}} \exp\left(-\frac{y^2}{4\alpha(t-\tau)}\right) d\tau \qquad [A11]$$

**Interfacial Heat Flux**

The interfacial heat flux at the fracture–rock boundary is

$$q(x,t) = -k \left.\frac{\partial T}{\partial y}\right|_{y=0} \qquad [A12]$$

Evaluating the derivatives yields

$$q(x,t) = \frac{k}{\sqrt{\pi\alpha}} \left[\frac{T_0}{\sqrt{t}} - \frac{1}{2} \int_0^t \frac{T_f(x,\tau)}{(t-\tau)^{3/2}} d\tau\right] \qquad [A13]$$

**Coupling with the Fluid Energy Equation**

Substitution of the interfacial heat flux into the fluid energy balance results in the integro-differential equation

$$\rho_f c_f v b \frac{\partial T_f}{\partial x} = -\frac{k}{\sqrt{\pi\alpha}} \left[\frac{T_0}{\sqrt{t}} - \frac{1}{2} \int_0^t \frac{T_f(x,\tau)}{(t-\tau)^{3/2}} d\tau\right] \qquad [A14]$$

**Solution in the Laplace Domain**

Introducing the Laplace transform $\tilde{T}_f(x,s) = \mathcal{L}\{T_f(x,t)\}$, the transformed equation becomes

$$\frac{\partial \tilde{T}_f}{\partial x} + \frac{k}{\rho_f c_f v b \sqrt{\alpha}} s^{1/2} \tilde{T}_f = -\frac{k T_0}{\rho_f c_f v b \sqrt{\alpha}} s^{-1/2} \qquad [A15]$$

This first-order linear ordinary differential equation admits the solution

$$\tilde{T}_f(x,s) = \frac{T_0}{s} + \frac{T_{\text{inj}} - T_0}{s} \exp\left(-\frac{k x}{\rho_f c_f v b \sqrt{\alpha}} \sqrt{s}\right) \qquad [A16]$$

where $T_{\text{inj}}$ is the fluid injection temperature.

**Inverse Laplace Transform**

Applying the inverse Laplace transform

$$\mathcal{L}^{-1}\left\{\frac{1}{s}\exp(a\sqrt{s})\right\} = \text{erfc}\left(\frac{a}{2\sqrt{t}}\right), \qquad a = \frac{kx}{\rho_f c_f v b \sqrt{\alpha}} \qquad [A17]$$

yields the fluid temperature distribution

$$T_f(x,t) = T_0 + (T_{\text{inj}} - T_0)\text{erfc}\left(\frac{a}{2\sqrt{t}}\right) \qquad [A18]$$

An equivalent representation is

$$T_f(x,t) = T_0 + (T_0 - T_{\text{inj}})\text{erf}\left(\frac{kx}{2\rho_f c_f v b \sqrt{\alpha t}}\right) \qquad [A19]$$

**Acknowledgements**

The author acknowledges support from the Utah FORGE project "High Temperature Proppants and Zeolite Markers: Designing, Characterizing & Optimizing Proppant and Flow Monitoring Materials for a Utah FORGE Engineered Geothermal System" (Project No. 9-3706). The insights gained through this project helped motivate the practical relevance of analytical models capable of informing geothermal system design and performance assessment.

The author also acknowledges the foundational contributions of Dr. A. C. Gringarten to analytical heat transfer modeling in fractured hot dry rock. The formulation presented in our work, primarily developed using Green's functions, is conceived as a scientific homage to the analytical tradition established by Dr. Gringarten, whose pioneering work laid the foundations for tractable analytical solutions to heat transfer problems in petroleum and geothermal reservoir engineering. While Green's functions provide a natural and elegant basis for the present formulation, the governing model can also be solved analytically using alternative methodologies, reinforcing its generality and robustness. His work has profoundly influenced multiple generations of researchers and practitioners.